\documentclass{article}

\usepackage{amsmath,yhmath,amsfonts,amssymb,amsthm,bbding}
\usepackage{mathrsfs}

\newcommand{\counte}{theorem}

\newtheorem{prop*}{\bf Step} 

\allowdisplaybreaks

\numberwithin{equation}{section}

\renewcommand{\thefootnote}{\fnsymbol{footnote}}

\begin{document}

\renewcommand{\thefootnote}{\arabic{footnote}}

\centerline{\bf\Large An Index I\!I\!I lemma and Rauch I\!I\!I theorem}
\centerline{\bf\Large\& applications \footnote{Supported by NSFC 11971057 and BNSF Z190003. \hfill{$\,$}}}

\vskip5mm

\centerline{Shengqi Hu, Xiaole Su,  Yusheng Wang\footnote{The
corresponding author (E-mail: wyusheng@bnu.edu.cn). \hfill{$\,$}}}

\vskip6mm

\noindent{\bf Abstract.} Inspired by Index I and I\!I lemmas and Rauch I and I\!I theorems, we formulate out an
Index I\!I\!I lemma and Rauch I\!I\!I theorem in this paper. As applications, we present a Rauch's type theorem with lower Ricci
curvature bound and a volume comparison result.

\vskip1mm


\vskip1mm


\vskip6mm

\setcounter{section}{0}


\section{Main results}

In Riemannian geometry, Rauch's type theorems are really fundamental comparison results essentially by sectional curvature ([Ra], [Be], [Wa]). In [CE], two of them
are listed in the sequel as Rauch I (Rauch) and Rauch I\!I (Berger) theorems, whose proofs need Index I and I\!I lemmas respectively. In [Wa], the Rauch theorems are generalized to
compare the Jacobi fields determined by smooth submanifolds with some assumptions on their second fundamental forms. Inspired by them, we will formulate out an
Index I\!I\!I lemma and Rauch I\!I\!I theorem in the paper.

Let $M$ be a Riemannian manifold, and let $\gamma:[0, l] \rightarrow M$ be a normal geodesic.
Given a symmetric linear transformation $B$ on $M_{\gamma(0)}$, if there is a Jacobi field
$J(t)$ along $\gamma$ with $J(l)=0$ but $0\neq J(0)\in M_{\gamma(0)}$ and $J'(0)=B(J(0))$, then we call $\gamma(l)$ a
{\it focal point of $\gamma(0)$ with respect to $B$ on $\gamma$}. Note that there is a $\lambda$ such that
$\frac{\|B(X)\|}{\|X\|}\leq\lambda$ for any nonzero $X\in M_{\gamma(0)}$,
so there is a $\delta\in (0, l]$ not depending on $X$ such that there is
a unique Jacobi field $X(t)$ along $\gamma|_{[0,\delta]}$ with $X(0)=X$, $X'(0)=B(X(0))$ and $X(t)\neq 0$,
i.e. there is no focal point of $\gamma(0)$ with respect to $B$ on $\gamma|_{[0,\delta]}$.

\vskip2mm

\noindent {\bf Lemma A.} (Index I\!I\!I). {\it Let $M$ be a Riemannian manifold, and let $\gamma:[0, l] \rightarrow M$ be a normal geodesic. Let $W$ be
a piecewise smooth vector field along $\gamma$, and $B$ be a symmetric linear transformation on $M_{\gamma(0)}$. If there is no focal point of $\gamma(0)$ with respect to $B$ on $\gamma$,
then there is a unique Jacobi field $V$ along $\gamma$ such that $V(l)=W(l)$ and $V^{\prime}(0)=B(V(0))$. Moreover,
$$I_B(V,V)\leq I_B(W,W),$$ and equality holds only if $V=W$.
}

\vskip2mm

In Lemma A,
$$I_B(Z,Z)\triangleq\langle B(Z(0)), Z(0)\rangle+\int_0^l\left(\left\langle Z', Z'\right\rangle+\langle R(\gamma',Z)\gamma', Z\rangle\right)dt.$$
If $Z$ is in addition a Jacobi field with $Z^{\prime}(0)=B(Z(0))$, then it is easy to see that
$$I_B(Z,Z)=\langle Z'(l), Z(l)\rangle.\eqno{(1.1)}$$

\vskip2mm

\noindent {\bf Remark 1.1.} A special case of Lemma A is that $B=\lambda\cdot\text{id}$, i.e. $V^{\prime}(0)=\lambda V(0)$; in particular, if $\lambda=0$, then the lemma is just the Index I\!I lemma in [CE].

\vskip2mm

\noindent {\bf Remark 1.2.} Lemma A should be known to experts. Note that if $V(0)$ and $W(0)$ lie in the tangent space of a submanifold $N$ of codimension 1 which is perpendicular to $\gamma$ at $\gamma(0)$, and if $V^{\prime}(0)=S_{\gamma'(0)}(V(0))$ (where $S_{\gamma'(0)}$ is the second fundamental form of $N$ with respect to $\gamma'(0)$ at $\gamma(0)$), then the lemma follows from the Minimization Theorem 2.1 in [Wa] (cf. [BC]).

\vskip2mm

As an application of Index I\!I\!I (Lemma A), we have a Rauch I\!I\!I theorem. Let $\gamma(t)|_{[0,l]}$ be a normal geodesic in a Riemannian manifold, and $V(t)$ be a smooth vector field along $\gamma$. We set $\widehat{V(t)}\triangleq V(t)-\langle V(t), \gamma'(t)\rangle \gamma'(t)$. Note that 
$(\widehat{V(t)})'=\widehat{V'(t)}$ for all $t$.

\vskip2mm

\noindent {\bf Theorem B.} (Rauch I\!I\!I). {\it Let $M, M_{0}$ be Riemannian manifolds, and let $\gamma, \gamma_{0}:[0, l] \rightarrow M, M_{0}$ be normal geodesics. Assume that for any $t \in[0, l]$, the sectional curvatures of sections spanned by $\gamma'(t), X\in M_{\gamma(t)}$ and $\gamma_0'(t), X_{0}\in\left(M_{0}\right)_{\gamma_{0}(t)}$ satisfy $K(\gamma'(t), X)\leq K\left(\gamma_0'(t),X_0\right)$. Let $V, V_{0}$ be Jacobi fields along $\gamma, \gamma_{0}$ satisfying $\|V(0)\|\geq\left\|V_{0}(0)\right\|>0$, $\langle \gamma'(0), V(0)\rangle=\left\langle \gamma'_{0}(0), V_{0}(0)\right\rangle$,  and $\langle \gamma'(0), V'(0)\rangle=\left\langle \gamma'_{0}(0), V_{0}'(0)\right\rangle$. If there are symmetric linear transformations $B, B_0$ on $\gamma'(0)^\perp,\gamma'_{0}(0)^\perp$ with the minimum eigenvalue of $B$ $\geq$ the maximum eigenvalue of $B_0$
such that $\widehat{V'(0)}=B(\widehat{V(0)}),\widehat{V_0^{\prime}(0)}=B_0(\widehat{V_0(0)})$, and if there is no focal point of $\gamma_0(0)$ with respect to $B_0$ on $\gamma_{0}$, then for all $t \in[0, l]$
$$\|V(t)\| \geq\left\|V_{0}(t)\right\|.\eqno{(1.2)}$$
Moreover, if $\|V(t_0)\|=\left\|V_{0}(t_0)\right\|$ for some $t_0\in(0,l]$, then, on $[0, t_0]$, $\|\widehat{V}\|=\|\widehat{V_0}\|$, $\|V\|=\left\|V_{0}\right\|$,
and $K(\gamma', V)=K\left(\gamma_0',V_0\right)$ and $\frac{\widehat{V}}{\|\widehat{V}\|},
\frac{\widehat{V_0}}{\|\widehat{V}\|}$ are parallel vector fields along $\gamma, \gamma_0$ if $\|\widehat{V(0)}\|\neq0$.
}

\vskip2mm

Note that in the case where $\|V(t_0)\|=\left\|V_{0}(t_0)\right\|$ in Theorem B, if $\widehat{V(0)}\neq 0$, then $\frac{\|\widehat{V(0)}\|'}{\|\widehat{V(0)}\|}$ is a common eigenvalue of $B$ and $B_0$,
and $\widehat{V(0)},\widehat{V_0(0)}$ are eigenvectors of $B,B_0$.

\vskip2mm

\noindent {\bf Remark 1.3.} Similar to Remark 1.1, a special case of Theorem B is that both $B$ and $B_0$ are $\lambda\cdot\text{id}$, i.e. $\widehat{V^{\prime}(0)}=\lambda \widehat{V(0)}$ and $\widehat{V_0^{\prime}(0)}=\lambda \widehat{V_0(0)}$; in particular, if $\lambda=0$, then the theorem is almost just the Rauch I\!I theorem in [CE] (but there it needs that $\dim(M)\leq\dim(M_0)$).

\vskip2mm

\noindent {\bf Remark 1.4.} Theorem B should also be known to experts except possibly some conclusions for `$\|V(t_0)\|=\left\|V_{0}(t_0)\right\|$'. Similarly, if $V(0),V_0(0)$ lie in the tangent spaces of submanifolds $N,N_0$ of codimension 1 which are perpendicular to $\gamma,\gamma_0$ at $\gamma(0),\gamma_0(0)$, and if $V^{\prime}(0)=S_{\gamma'(0)}(V(0))$, $V_0^{\prime}(0)=S_{\gamma_0'(0)}(V_0(0))$, and the minimum eigenvalue of $S_{\gamma'(0)}$ is not less than the maximum eigenvalue of $S_{\gamma_0'(0)}$, then (1.2) follows from the case (b) of Theorem 4.3 in [Wa].

\vskip2mm

\noindent {\bf Remark 1.5.} We would like to point out that, about the case `$\|V(t_0)\|=\left\|V_{0}(t_0)\right\|$',
Theorem 4.3 in [Wa] just say that $\|V(t)\|=\left\|V_{0}(t)\right\|$ on $[0, t_0]$ (and Rauch I and I\!I theorems in [CE] say nothing), which follows immediately from the proof for `$\|V(t)\|\geq\left\|V_{0}(t)\right\|$'.
That we can see `$\frac{\widehat{V(t)}}{\|\widehat{V(t)}\|}$
and $\frac{\widehat{V_0(t)}}{\|\widehat{V(t)}\|}$ are both parallel vector fields' in Theorem B is due to a quite different trick from [CE] and [Wa] (which can yield the same rigidity for Rauch I). Moreover, the trick also guarantees that our proof goes through whether or not $\dim(M)\leq\dim(M_0)$ (but it is needed in
Rauch I and I\!I).

\vskip2mm

In [CE], two interesting corollaries of Rauch I and I\!I theorems are listed (cf. [Be]). Similarly, as an application of the special case where $V^{\prime}(0)=\lambda V(0)$ and $V_0^{\prime}(0)=\lambda V_0(0)$ in Theorem B,
we  have the following corollary.

\vskip2mm

\noindent {\bf Corollary C.} {\it Let $M, M_{0}$ be Riemannian manifolds with sectional curvature $K_{M}\leq K_{M_{0}}$, and let $\gamma, \gamma_{0}:[0, l] \rightarrow M, M_{0}$ be normal geodesics.
Let $E,E_0$ be smooth vector fields along $\gamma,\gamma_{0}$ such that
$\|E(t)\|=\|E_0(t)\|\neq 0$, $\langle E(t), \gamma'(t)\rangle=\left\langle E_{0}(t), \gamma'_{0}(t)\right\rangle\neq\pm\|E(t)\|$, $E'(t)=\lambda(t)\gamma'(t)$ and $E_0'(t)=\lambda(t)\gamma'_0(t)$
on $[0,l]$. Let $c, c_0:[0, l] \rightarrow M, M_{0}$ be  smooth curves defined by
\[ \\
c(t)=\exp_{\gamma(t)}(f(t) E(t)) \text{ and } c_{0}(t)=\exp_{\gamma_0(t)}\left(f(t) E_{0}(t)\right),\\
\]
where $f$ is a non-negative smooth function. Then if, for any fixed $t$, there is no focal point of $\gamma_0(t)$  with respect to $\lambda(t)\cdot \text{\rm id}$ on $\exp_{\gamma_0(t)}\left(s E_{0}(t)/\|E_{0}(t)\|\right)|_{[0,f(t)\|E_{0}(t)\|]}$, then the lengths of $c$ and $c_0$ satisfy
\[ \\
L[c] \geq L\left[c_{0}\right].
\]}
\hskip2.5mm In Corollary C, if in addition $\lambda(t)\equiv0$ and $E,E_{0}$ are perpendicular to $\gamma',\gamma'_{0}$, then the corollary is just the Corollary of Rauch I\!I in [CE]
(cf. Theorem 1 in [Be]).

\vskip2mm

\noindent {\bf Remark 1.6.} As an easy example of Corollary C, we consider minimal geodesics $[pq], [pr], [qs]\subset \Bbb S^2_0$ (a plane) and $[p_0q_0],[p_0r_0], [q_0s_0]\subset \Bbb S^2_1$ (a unit sphere)
with the distance $|pq|=|p_0q_0|<\pi$, $|pr|=|p_0r_0|$ and $|qs|=|q_0s_0|$. Assume that the angle $\angle rpq=\angle r_0p_0q_0$, $\angle pqs=\angle p_0q_0s_0$, and $r$ and $s$ (resp. $r_0$ and $s_0$) lie in one side of $[pq]$
(resp. $[p_0q_0]$). Then there is $C_1>0$ and $C_2>0$ depending on $|pq|$, $\angle rpq$ and $\angle pqs$  such that if $|pr|\leq C_1$ and $|qs|\leq C_2$ then
$$|rs|\geq |r_0s_0|.$$
(In fact, the $\Bbb S^2_0, \Bbb S^2_1$ in the example can be replaced with 2-dimensional complete Riemannian manifolds $M, M_{0}$ with Gauss curvature $K_{M}\leq K_{M_{0}}$.)

\vskip2mm

As another application of Lemma A, we can get a Rauch's type theorem with lower Ricci curvature bound which just relies on the special case where $V^{\prime}(0)=\lambda V(0)$ in the lemma.

\vskip2mm

\noindent {\bf Theorem D.} {\it Let $M^{n}$ be a Riemannian manifold with $\text{\rm Ric}_M\geq (n-1)k$, and let $\Bbb S^n_k$ be the simply connected and complete space form of constant curvature $k$. Let $\gamma, \tilde\gamma:[0, l] \rightarrow M^n, \Bbb S^n_k$ be normal geodesics. Let $V_i, \tilde V_{i}$, $1\leq i\leq n-1$, be Jacobi fields along $\gamma, \tilde \gamma$ with $\langle \gamma', V_i\rangle=\langle \tilde \gamma', \tilde V_{i}\rangle=0$ and $0<\|V_1(0)\wedge\cdots\wedge V_{n-1}(0)\|\leq\|\tilde V_{1}(0)\wedge\cdots\wedge\tilde V_{n-1}(0)\|$. If $V_i^{\prime}(0)=\lambda V_i(0)$ and $\tilde V_{i}^{\prime}(0)=\tilde\lambda \tilde V_{i}(0)$ with $\lambda\leq\tilde\lambda$, and if there is no focal point of $\gamma(0)$ with respect to $\lambda\cdot\text{\rm id}$ on $\gamma$, then for all $t \in[0, l]$
$$\|\tilde V_1(t)\wedge\cdots\wedge \tilde V_{n-1}(t)\| \geq\left\|V_{1}(t)\wedge\cdots\wedge V_{n-1}(t)\right\|;$$
moreover, if the equality holds for some $t_0\in(0,l]$, then on $[0, t_0]$ each $\|V_i\|=\|\tilde V_{i}\|$, $\frac{V_{i}}{\|V_{i}\|}$ is a parallel vector field,
$K(\gamma',V_i)=k$, and $\lambda=\tilde\lambda$.}

\vskip2mm

Note that, in Theorem D, $\frac{\tilde V_{i}}{\|\tilde V_{i}\|}$ is naturally a parallel vector field. And it turns out that Theorem D can be reduced to `$\{\gamma'(0),V_1(0),\cdots,V_{n-1}(0)\},\{\tilde \gamma'(0),\tilde V_{1}(0),\cdots,\tilde V_{n-1}(0)\}$ are orthonormal bases of $M_{\gamma(0)},(\Bbb S^n_k)_{\tilde\gamma(0)}$', a special case of `$0<\|V_1(0)\wedge\cdots\wedge V_{n-1}(0)\|\leq\|\tilde V_{1}(0)\wedge\cdots\wedge\tilde V_{n-1}(0)\|$'.

\vskip2mm

\noindent {\bf Remark 1.7.} In Theorem D, if the initial values of $V_i,\tilde V_i$ ($1\leq i\leq n-1$) are $V_i(0)=0$, $\tilde V_i(0)=0$, and
$|\gamma'(0)\wedge V_1'(0)\wedge\cdots\wedge V_{n-1}'(0)|=|\tilde\gamma'(0)\wedge\tilde V_{1}'(0)\wedge\cdots\wedge\tilde V_{n-1}'(0)|$, and if there is no conjugate point of $\gamma(0)$ on $\gamma$,
then it also holds that $|\tilde\gamma'(t)\wedge\tilde V_1(t)\wedge\cdots\wedge \tilde V_{n-1}(t)|\geq\left|\gamma'(t)\wedge V_{1}(t)\wedge\cdots\wedge V_{n-1}(t)\right|$ (cf. [WSY]).
This can be viewed as a generalization of Rauch I theorem, while Theorem D can be viewed as a generalization of Rauch  I\!I and  I\!I\!I.
Moreover, in any case, we in fact have that $\frac{\left|\gamma'(t)\wedge V_{1}(t)\wedge\cdots\wedge V_{n-1}(t)\right|}{|\tilde\gamma'(t)\wedge\tilde V_1(t)\wedge\cdots\wedge \tilde V_{n-1}(t)|}$ is decreasing with respect to $t$.

\vskip2mm

As a corollary of Theorem D, we have a volume comparison result.

\vskip2mm

\noindent {\bf Corollary E.} {\it Let $M^{n}$ be a Riemannian manifold with $\text{\rm Ric}_M\geq (n-1)k$, and let $p\in M$ and $\tilde p\in \Bbb S^n_k$.
Suppose that $\overline{B(p,r)}$ is isometric to $\overline{B(p',r)}\subset \Bbb S^n_{k'}$ with $r>0$ and $k'\geq k$, and that $\text{\rm Vol}(\partial B(p,r))=\text{\rm Vol}(\partial B(\tilde p,\tilde r))$ with $\tilde r\leq\frac{\pi}{2\sqrt k}$ if $k>0$.
Then for any $R>0$,
$$\text{\rm Vol}(\partial B(p,r+R))\leq\text{\rm Vol}(\partial B(\tilde p,\tilde r+R)),$$
which implies $\text{\rm Vol}\left(B(p,r+R)\setminus B(p,r)\right)\leq\text{\rm Vol}\left(B(\tilde p,\tilde r+R)\setminus B(\tilde p,\tilde r)\right)$ and the equality holds only if $B(p,r+R)$ is
isometric to $B(\tilde p,\tilde r+R)$. }

\vskip2mm

In Corollary E, if $r=0$, then it is just the classical Bishop's volume comparison theorem ([Pe], [WSY]), which can be deduced not hardly from the result in Remark 1.7
\footnote{The Bishop-Gromov relative volume comparison theorem relies on the monotonicity of $\frac{\left|\gamma'(t)\wedge J_{1}(t)\wedge\cdots\wedge J_{n-1}(t)\right|}{|\tilde\gamma'(t)\wedge\tilde J_1(t)\wedge\cdots\wedge \tilde J_{n-1}(t)|}$ in Remark 1.7 ([Pe], [WSY]).}.

\section{Proofs of main results}

As mentioned above, Lemma A should be known to experts. Note that for any $X\in M_{\gamma(0)}$ in Lemma A, there is a unique Jacobi field $X(t)$ on $[0,l]$ with $X(0)=X$ and $X'(0)=B(X(0))$,
and the map $f:M_{\gamma(0)}\to M_{\gamma(l)}$ defined by $f(X)=X(l)$ is linear.
Since there is no focal point of $\gamma(0)$ with respect to $B$ on $\gamma$, $\ker(f)=\{0\}$ and thus $f$ has to
be a 1-1 map. As a result, given a basis $\{V_i\}$ of $M_{\gamma(l)}$,
we can extend each $V_i$ to a Jacobi field $V_i(t)$ on $[0,l]$ with $V_i'(0)=B(V_i(0))$ and $\{V_i(t)\}$ being linearly independent.
And thus if we write $W(l)=\sum c_iV_i(l)$, then $V(t)\triangleq\sum c_iV_i(t)$ on $[0,l]$ is the unique Jacobi field satisfying
$V(l)=W(l)$ and $V^{\prime}(0)=B(V(0))$. In such a situation, the rest of the proof of Lemma A is almost a copy of that of Index I lemma in [CE].

\vskip2mm

As for Theorem B, although (1.2) might be known to experts, we would like to provide a detailed proof for the theorem, in which a different trick from [CE] and [Wa] enables us to see `$\frac{\widehat{V(t)}}{\|\widehat{V(t)}\|}$
and $\frac{\widehat{V_0(t)}}{\|\widehat{V(t)}\|}$ are parallel vector fields' in the theorem.
And the trick guarantees that our arguments have nothing to do with the dimensions of $M$ and $M_0$.

\vskip2mm

\noindent {\bf Proof of Theorem B}.

We first consider the case where  $\langle \gamma'(0), V(0)\rangle=\langle \gamma'(0), V'(0)\rangle=0$. In this case,
$V(t)=\widehat{V(t)}$ (i.e. $\langle \gamma'(t), V(t)\rangle\equiv0$) and $V_0(t)=\widehat{V_0(t)}$.

Note that $V_0(t)\neq 0$ on $[0,l]$ because $V_0(0)\neq 0$, $V_0^{\prime}(0)=B_0(V_0(0))$ and there is no focal point of $\gamma_0(0)$ with respect to $B_0$ on $\gamma_{0}$. So, as in the proof of Rauch I theorem in [CE], we can consider the ratio function \(\frac{\|V(t)\|^{2}}{\left\|V_{0}(t)\right\|^{2}}\) on \([0,l]\) with $\frac{\|V(0)\|^{2}}{\left\|V_{0}(0)\right\|^{2}}\geq 1$.
Then it is not hard to see that, to prove \(\|V(t)\| \geq\left\|V_{0}(t)\right\|\) on $[0,l]$, it suffices to show that, for \(t_1\in (0,l)\) with $V(t)\neq 0$ on $[0,t_1]$ \footnote{In the proof of Theorem 4.3 in [Wa], it is first proven that $V(t)\neq 0$ on $(0,l]$, while this is just a natural corollary in our deduction.},
\[ \\
\left.\frac{\mathrm{d}}{\mathrm{d} t}\left(\frac{\|V\|^{2}}{\left\|V_{0}\right\|^{2}}\right)\right|_{t_1} \geq 0, \text{ or equivalently, }
\left.\frac{\left\langle V^{\prime}, V\right\rangle}{\langle V, V\rangle}\right|_{t_1}  \geq \left.\frac{\left\langle V_{0}^{\prime}, V_{0}\right\rangle}{\left\langle V_{0}, V_{0}\right\rangle}\right|_{t_1}. \eqno{(2.1)}
\]
For the purpose, we consider Jacobi fields on $[0,t_1]$ (as in [CE])
\[ \\
W_{t_{1}}(t)=\frac{V(t)}{\left\|V\left(t_{1}\right)\right\|},\ \ W_{0 t_{1}}(t)=\frac{V_{0}(t)}{\left\|V_{0}\left(t_{1}\right)\right\|} . \\
\] \\
Then \(\left\|W_{t_1}\left(t_{1}\right)\right\|=\left\|W_{0t_1}\left(t_{1}\right)\right\|=1\), and
\[ \\
\frac{\left\langle V', V\right\rangle}{\langle V, V\rangle}=\frac{\langle W_{t_{1}}', W_{t_{1}}\rangle}{\left\langle W_{t_{1}}, W_{t_{1}}\right\rangle},\ \ \frac{\left\langle V_{0}', V_{0}\right\rangle}{\left\langle V_{0}, V_{0}\right\rangle}=\frac{\langle W_{0 t_{1}}', W_{0 t_{1}}\rangle}{\left\langle W_{0 t_{1}}, W_{0 t_{1}}\right\rangle} . \\
\] \\
In particular, \\
\[ \\
\left.\frac{\left\langle V^{\prime}, V\right\rangle}{\langle V, V\rangle}\right|_{t_{1}}=\left.\left\langle W^{\prime}_{t_{1}}, W_{t_{1}}\right\rangle\right|_{t_{1}},\left.\ \ \frac{\left\langle V^{\prime}_{0}, V_{0}\right\rangle}{\left\langle V_{0}, V_{0}\right\rangle}\right|_{t_{1}}=\left.\left\langle W^{\prime}_{0 t_{1}}, W_{0 t_{1}}\right\rangle\right|_{t_{1}}. \\
\] \\
Note that $\left.\left\langle W_{t_{1}}^{\prime}, W_{t_{1}}\right\rangle\right|_{t_{1}}=I_B(W_{t_{1}},W_{t_{1}})$
and $\left.\left\langle W_{0t_{1}}^{\prime}, W_{0t_{1}}\right\rangle\right|_{t_{1}}
=I_{B_0}(W_{0t_{1}},W_{0t_{1}})$ (see (1.1)). I.e, we need to show that $$I_B(W_{t_{1}},W_{t_{1}})\geq I_{B_0}(W_{0t_{1}},W_{0t_{1}}).\eqno{(2.2)}$$

Let $e_0(t)$, $t\in[0,t_1]$, be the parallel (unit) vector field along $\gamma_0$ with $e_0(t_1)=W_{0t_{1}}(t_1)$; and let $\overline{W}_0(t)\triangleq\|W_{t_{1}}(t)\|e_0(t)$
\footnote{This construction is just the trick mentioned right before the Proof of Theorem B.}.
Note that $\overline{W}_0(t_1)=W_{0t_{1}}(t_1)$. Then, by Lemma A we have that
$$I_{B_0}(\overline{W}_0,\overline{W}_0)\geq I_{B_0}(W_{0t_{1}},W_{0t_{1}}).$$
On the other hand, using our assumption on the curvatures and eigenvalues gives \\
\[
\begin{aligned}
I_{B_0}(\overline{W}_0,\overline{W}_0)
&=\left\langle B_0\left(\overline{W}_0(0)\right), \overline{W}_0(0)\right\rangle+
  \int_{0}^{t_{1}}\left(\left\langle \overline{W}_0^{\prime}, \overline{W}_0^{\prime}\right\rangle+\left\langle R\left(\gamma'_0, \overline{W}_0\right) \gamma'_0, \overline{W}_0\right\rangle\right)dt \\
& \leq \left\langle B\left(W_{t_1}(0)\right), W_{t_1}(0)\right\rangle+
   \int_{0}^{t_{1}}\left(\left(\|W_{t_{1}}(t)\|'\right)^2+\left\langle R\left(\gamma', W_{t_1}\right) \gamma', W_{t_1}\right\rangle\right)dt \\
& \leq \left\langle B\left(W_{t_1}(0)\right), W_{t_1}(0)\right\rangle+
   \int_{0}^{t_{1}}\left(\left\langle W_{t_{1}}^{\prime}, W_{t_{1}}^{\prime}\right\rangle+\left\langle R\left(\gamma', W_{t_1}\right) \gamma', W_{t_1}\right\rangle\right)dt. \\
& =I_B(W_{t_{1}},W_{t_{1}}).
\end{aligned}
\]
(Here, for the first `$\leq$', note that $\|\overline{W}_0(t)\|=\|W_{t_1}(t)\|$, and $\overline{W}_0,W_{t_1}$ are perpendicular to $\gamma_0,\gamma$;
for the second one, note that $e(t)\triangleq\frac{W_{t_1}(t)}{\|W_{t_1}(t)\|}$ is a smooth unit vector field, and thus $\left\langle W_{t_{1}}^{\prime}, W_{t_{1}}^{\prime}\right\rangle=\left(\|W_{t_{1}}(t)\|'\right)^2+\|W_{t_{1}}(t)\|^2\|e'(t)\|^2$.)
It then follows that $$I_B(W_{t_{1}},W_{t_{1}})\geq I_{B_0}(\overline{W}_0,\overline{W}_0)\geq I_{B_0}(W_{0t_{1}},W_{0t_{1}}),\eqno{(2.3)}$$
i.e. (2.2) and so (2.1) hold. Moreover, the equality in (2.1), or equivalently the two equalities in (2.3), implies that   $W_{0t_{1}}(t)=\overline{W}_0(t)$, $e'(t)=0$,
and $\left\langle R\left(\gamma'_0, \overline{W}_0\right) \gamma'_0, \overline{W}_0\right\rangle=\left\langle R\left(\gamma', W_{t_1}\right) \gamma', W_{t_1}\right\rangle$  on $[0,t_1]$.

So far, we have proven that \(\|V\| \geq\left\|V_{0}\right\|>0\) on $[0,l]$ with $\left(\frac{\|V\|^{2}}{\left\|V_{0}\right\|^{2}}\right)'\geq 0$ (see (2.1)). Thereby, if $\|V(t_0)\|=\left\|V_{0}(t_0)\right\|$ for some $t_0\in(0,l]$, then $\|V\|=\|V_0\|$ and $\left(\frac{\|V\|^{2}}{\left\|V_{0}\right\|^{2}}\right)'=0$ on $[0, t_0]$; and it follows that both $\frac{V}{\|V\|}$
and $\frac{V_0}{\|V\|}$ are parallel vector fields, and $K(\gamma',V)=K\left(\gamma_0',V_0\right)$ on $[0, t_0]$.

\vskip2mm

We now consider general cases where $V(t)=\widehat{V(t)}+(at+b)\gamma'(t)$ and $V_{0}(t)=\widehat{V_{0}(t)}+(at+b)\gamma'_0(t)$ with $b=\langle \gamma'(0), V(0)\rangle$ and $a=\langle \gamma'(0), V'(0)\rangle$. Note that if $\widehat{V_0(0)}=0$ and so $(\widehat{V_0(0)})'=B(\widehat{V_0(0)})=0$, then the Jacobi field $\widehat{V_0(t)}\equiv0$, so does $\widehat{V(t)}$; and thus
$\|V(t)\|=\left\|V_{0}(t)\right\|$ on $[0,l]$. If $\widehat{V_0(0)}\neq 0$, then we can apply the above arguments to $\widehat{V(t)}$ and $\widehat{V_0(t)}$, and complete the proof.
\hfill$\Box$

\vskip2mm

\noindent {\bf Proof of Corollary C}.

It suffices to show that, for any $t_0\in(0,l)$ with $f(t_0)>0$,
$$\|c'(t_0)\|\geq \|c_0'(t_0)\|.$$
Since $\|E(t_0)\|=\|E_0(t_0)\|\neq0$, for sufficiently small $\delta>0$ we can define a smooth map
$$\sigma:(t_0-\delta,t_0+\delta)\times[0,f(t_0)\|E(t_0)\|]\to M \text{ by } \sigma(t,s)=\exp_{\gamma(t)}\left(\frac{sf(t)}{f(t_0)}\frac{E(t)}{\|E(t_0)\|}\right),$$
and similarly $\sigma_0(t,s)=\exp_{\gamma_0(t)}\left(\frac{sf(t)}{f(t_0)}\frac{E_0(t)}{\|E(t_0)\|}\right)$. It is clear that $\sigma(t,f(t_0)\|E(t_0)\|)=c(t)|_{(t_0-\delta,t_0+\delta)}$
and $\|c'(t_0)\|=\left\|\frac{\partial}{\partial t}\sigma(t_0,f(t_0)\|E(t_0)\|)\right\|$, and similarly for $\sigma_0(t,f(t_0)\|E(t_0)\|)$.
Since $\sigma$ and $\sigma_0$ are both smooth, $\left.\nabla_{\frac{\partial}{\partial t}}\frac{\partial}{\partial s}\right|_{(t_0,0)}=\left.\nabla_{\frac{\partial}{\partial s}}\frac{\partial}{\partial t}\right|_{(t_0,0)}$ on $M$ and $\left.{\nabla_0}_{\frac{\partial}{\partial t}}\frac{\partial}{\partial s}\right|_{(t_0,0)}=\left.{\nabla_0}_{\frac{\partial}{\partial s}}\frac{\partial}{\partial t}\right|_{(t_0,0)}$ on $M_0$ with
$$\begin{aligned}\left.\nabla_{\frac{\partial}{\partial t}}\frac{\partial}{\partial s}\right|_{(t_0,0)}&
=\left.\left(\frac{f(t)}{f(t_0)}\frac{E(t)}{\|E(t_0)\|}\right)'\right|_{t_0}=\frac{f'(t_0)}{f(t_0)}\frac{E(t_0)}{\|E(t_0)\|}
 +\frac{\lambda(t_0)}{\|E(t_0)\|}\gamma'(t_0),\\
\left.{\nabla_0}_{\frac{\partial}{\partial t}}\frac{\partial}{\partial s}\right|_{(t_0,0)}&=\frac{f'(t_0)}{f(t_0)}\frac{E_0(t_0)}{\|E(t_0)\|}+\frac{\lambda(t_0)}{\|E(t_0)\|}\gamma'_0(t_0).
\end{aligned}$$
Therefore, the Jacobi fields $V(s), V_0(s)\triangleq\left.\frac{\partial}{\partial t}\right|_{\{t_0\}\times[0,f(t_0)\|E(t_0)\|]}$ along $\sigma(t_0,s),\sigma_0(t_0,s)$ satisfy
$$V'(0)=\frac{f'(t_0)}{f(t_0)}\frac{E(t_0)}{\|E(t_0)\|}
 +\frac{\lambda(t_0)}{\|E(t_0)\|}\gamma'(t_0),\ \ V_0'(0)=\frac{f'(t_0)}{f(t_0)}\frac{E_0(t_0)}{\|E(t_0)\|}+\frac{\lambda(t_0)}{\|E(t_0)\|}\gamma'_0(t_0).$$
This plus $\langle E(t_0), \gamma'(t_0)\rangle=\langle E_0(t_0), \gamma'_0(t_0)\rangle\neq\pm\|E(t_0)\|$ implies that there is an $a$ and $b$ such that $V(s)=\widehat{V(s)}+(as+b)\frac{\partial}{\partial s}|_{(t_0,s)}, V_0(s)=\widehat{V_0(s)}+(as+b)\frac{\partial}{\partial s}|_{(t_0,s)}$ with $\widehat{V(0)}\neq0,\widehat{V_0(0)}\neq 0$
(note that $V(0)=\gamma'(t_0), V_0(0)=\gamma_0'(t_0)$), and thus
$$(\widehat{V(0)})'=\left(\frac{f'(t_0)}{f(t_0)}-a\right)\frac{E(t_0)}{\|E(t_0)\|}+\frac{\lambda(t_0)}{\|E(t_0)\|}V(0)$$
(similarly for $(\widehat{V_0(0)})'$).
On the other hand, note that $(\widehat{V(0)})'=\widehat{V'(0)}$ (see the comments right above Theorem B), so $(\widehat{V(0)})'$ is parallel to $\widehat{V(0)}$ (similarly for $(\widehat{V_0(0)})'$). Then it has to hold that
$$(\widehat{V(0)})'=\frac{\lambda(t_0)}{\|E(t_0)\|}\widehat{V(0)},\ \ (\widehat{V_0(0)})'=\frac{\lambda(t_0)}{\|E(t_0)\|}\widehat{V_0(0)}.$$
Since there is no focal point of $\gamma_0(t_0)$ with respect to $\lambda(t_0)\cdot\text{id}$ on $\sigma_0(t_0,s)|_{[0,f(t_0)\|E(t_0)\|]}$, by Theorem B \footnote{Here, it will be easier to apply Theorem B than Theorem 4.3 in [Wa].} we can conclude that
$$\|V(f(t_0)\|E(t_0)\|)\|\geq \|V_0(f(t_0)\|E(t_0)\|)\|,$$
i.e. $\|c'(t_0)\|\geq \|c_0'(t_0)\|$.
\hfill$\Box$

\vskip2mm

\noindent {\bf Proof of Theorem D}.

By Theorem B, we need only consider the case where $n\geq 3$. Let $J_i, \tilde J_{i}$, $1\leq i\leq n-1$, be Jacobi fields along $\gamma, \tilde \gamma$ such that $\langle \gamma', J_i\rangle=\langle \tilde \gamma', \tilde J_{i}\rangle=0$,
$\{\gamma'(0),J_1(0),$ $\cdots,J_{n-1}(0)\},\{\tilde \gamma'(0),\tilde J_{1}(0),\cdots,\tilde J_{n-1}(0)\}$ are orthonormal bases of $M_{\gamma(0)},(\Bbb S^n_k)_{\tilde\gamma(0)}$ \footnote{If $J_i, \tilde J_{i}$ satisfy only $\langle \gamma'(0), J_i(0)\rangle=\langle \tilde \gamma'(0), \tilde J_{i}(0)\rangle$ and $\langle \gamma'(0), J'_i(0)\rangle=\langle \tilde \gamma'(0), \tilde J'_{i}(0)\rangle$ instead of $\langle \gamma', J_i\rangle=\langle \tilde \gamma', \tilde J_{i}\rangle=0$, it needs that $\{\gamma'(0),\widehat{J_1(0)},\cdots\},\{\tilde \gamma'(0),\widehat{\tilde J_{1}(0)},\cdots\}$ are orthonormal bases of $M_{\gamma(0)},(\Bbb S^n_k)_{\tilde\gamma(0)}$.},
and $J_i^{\prime}(0)=\lambda J_i(0),\tilde J_{i}^{\prime}(0)=\tilde\lambda \tilde J_{i}(0)$.
Since there is no focal point of $\gamma(0)$ with respect to $\lambda\cdot\text{id}$ on $\gamma$, $\{\gamma'(t),J_1(t),\cdots,J_{n-1}(t)\}$ are also bases of $M_{\gamma(t)}$
on $[0,l]$, so Jacobi fields $V_1,\cdots,V_{n-1}$ can be linearly represented by $J_1,\cdots,J_{n-1}$.
Similarly, if $\{\tilde\gamma'(t),\tilde J_1(t),\cdots,\tilde J_{n-1}(t)\}$ are bases on $[0,l]$, then
$\tilde V_1,\cdots,\tilde V_{n-1}$ can be also linearly represented by $\tilde J_1,\cdots,\tilde J_{n-1}$.
It follows that
$$\begin{aligned}\|V_1(t)\wedge\cdots\wedge V_{n-1}(t)\|&=\|V_1(0)\wedge\cdots\wedge V_{n-1}(0)\|\cdot\|J_1(t)\wedge\cdots\wedge J_{n-1}(t)\|,\\
\|\tilde V_1(t)\wedge\cdots\wedge \tilde V_{n-1}(t)\|&=\|\tilde V_1(0)\wedge\cdots\wedge \tilde V_{n-1}(0)\|\cdot\|\tilde J_1(t)\wedge\cdots\wedge \tilde J_{n-1}(t)\|.
\end{aligned}$$
Therefore, it suffices to show that, for all $t \in[0, l]$,
$$\|\tilde J_1(t)\wedge\cdots\wedge \tilde J_{n-1}(t)\| \geq\left\|J_{1}(t)\wedge\cdots\wedge J_{n-1}(t)\right\|;$$ moreover, if the equality holds for $t_0\in(0,l]$, then on $[0, t_0]$ each $\|J_i\|=\|\tilde J_{i}\|$, $\frac{J_{i}}{\|J_{i}\|}$ is a parallel vector field \footnote{Note that each $\frac{\tilde J_{i}}{\|\tilde J_{i}\|}$ is naturally a parallel vector field and $\|\tilde J_1\|=\cdots=\|\tilde J_{n-1}\|$.}
(which together with $\|V_1(0)\wedge\cdots\wedge V_{n-1}(0)\|>0$ implies that each $\frac{V_{i}}{\|V_{i}\|}$ is a parallel vector field), $K(\gamma',J_i)=k$, and $\lambda=\tilde\lambda$.

For convenience, we let $J(t)\triangleq\left\|J_{1}(t)\wedge\cdots\wedge J_{n-1}(t)\right\|$ and $\tilde J(t)\triangleq\|\tilde J_1(t)\wedge\cdots\wedge \tilde J_{n-1}(t)\|$.
Note that $J(0)=\tilde J(0)\neq 0$ and  $J(t)\neq 0$ on $[0,l]$. Then it is not hard to see that, to prove $\tilde J(t)\geq J(t)$ on $[0,l]$, it suffices to show that, for \(t_1\in (0,l)\) with $\tilde J(t)\neq 0$ on $[0,t_1]$,
\[ \\
\left.\frac{\mathrm{d}}{\mathrm{d} t}\left(\frac{\tilde J(t)}{J(t)}\right)\right|_{t_1} \geq 0, \text{ or equivalently, }
\frac{\tilde J'(t_1)}{\tilde J(t_1)} \geq \frac{J'(t_1)}{J(t_1)}.\eqno{(2.4)}
\]

Let $\{e_i(t)\}_{i=1}^{n-1}$ be unit parallel vector fields along $\gamma$ such that $e_i(0)=J_i(0)$ and $J_i(t)=\sum\limits_{j=1}^{n-1}a_{ij}(t)e_j(t)$ with $J(t)=|(a_{ij}(t))|$ (the determinant). Since there is no focal point of $\gamma(0)$ with respect to $\lambda\cdot\text{id}$ on $\gamma$, there is unique Jacobi field $\hat e_i(t)$ along $\gamma|_{[0,t_1]}$ such that $\hat e_i(t_1)=e_i(t_1)$ and $\hat e_i'(0)=\lambda\hat e_i(0)$. {\bf Claim}:
$$\frac{J'(t_1)}{J(t_1)}=\sum_{i=1}^{n-1}I_{\lambda\cdot\text{id}}(\hat e_i,\hat e_i).\eqno{(2.5)}$$
In fact, $e_i(t)=\sum\limits_{j=1}^{n-1}a^{ij}(t)J_j(t)$ with $(a^{ij}(t))=(a_{ij}(t))^{-1}$, and thus the Jacobi field $\hat e_i(t)=\sum\limits_{j=1}^{n-1}a^{ij}(t_1)J_j(t)$. Then by a straight computation, we have
$$\begin{aligned}\frac{J'(t_1)}{J(t_1)}=& \frac{1}{J(t_1)}\sum_{i=1}^{n-1}
  \left|\begin{matrix} a_{11}(t_1) & \cdots& a_{1i}'(t_1)&\cdots&a_{1,n-1}(t_1)\\
                       \cdots & \cdots& \cdots&\cdots&\cdots\\
                       a_{n-1,1}(t_1) & \cdots& a_{n-1,i}'(t_1)&\cdots&a_{n-1,n-1}(t_1)\\
  \end{matrix}\right|\\
  =& \sum_{i=1}^{n-1}\sum_{j=1}^{n-1} a^{ij}(t_1)a_{ji}'(t_1)\\
  =& \sum_{i=1}^{n-1}\sum_{j=1}^{n-1} a^{ij}(t_1)\left.\left\langle J_j(t), e_i(t)\right\rangle'\right|_{t_1}\\
  =& \sum_{i=1}^{n-1}\sum_{j=1}^{n-1} \left\langle a^{ij}(t_1)J_j'(t_1), e_i(t_1)\right\rangle\\
  =& \sum_{i=1}^{n-1}\left\langle \hat e_i'(t_1), \hat e_i(t_1)\right\rangle.
 \end{aligned} $$
Note that $\langle\hat e_i'(t_1), \hat e_i(t_1)\rangle=I_{\lambda\cdot\text{id}}(\hat e_i,\hat e_i)$ (see (1.1)) because $\hat e_i(t)$ is a Jacobi field with $\hat e_i'(0)=\lambda\hat e_i(0)$, and so the claim follows.

On the other hand, on $\Bbb S^n_k$, we have that $\|\tilde J_1(t)\|=\cdots=\|\tilde J_{n-1}(t)\|\triangleq a(t)$ on $[0,t_1]$, and each $\tilde e_i(t)\triangleq\frac{\tilde J_i(t)}{a(t)}$ is a parallel vector field, and $\hat{\tilde e}_i(t)\triangleq\frac{\tilde J_i(t)}{a(t_1)}=\frac{a(t)}{a(t_1)}\tilde e_i(t)$ is also a Jacobi field. As a result, $\tilde J(t)=a^{n-1}(t)$, and so
$$\frac{\tilde J'(t_1)}{\tilde J(t_1)}=(n-1)\frac{a'(t_1)}{a(t_1)}
=\sum_{i=1}^{n-1}\langle \hat{\tilde e}_i'(t_1), \hat{\tilde e}_i(t_1)\rangle
=\sum_{i=1}^{n-1}I_{\tilde\lambda\cdot\text{id}}(\hat{\tilde e}_i,\hat{\tilde e}_i). \eqno{(2.6)}$$

We now set $\bar e_i(t)\triangleq \frac{a(t)}{a(t_1)} e_i(t)$ along $\gamma(t)|_{[0,t_1]}$, $i=1,\cdots, n-1$.
Note that $\bar e_i(t_1)=\hat e_i(t_1)$. Then by Lemma A and our assumption on the curvatures, we have that
$$
\begin{aligned}
\sum_{i=1}^{n-1}I_{\lambda\cdot\text{id}}(\hat e_i,\hat e_i)
&\leq \sum_{i=1}^{n-1}I_{\lambda\cdot\text{id}}(\bar e_i,\bar e_i) \\
& = \sum_{i=1}^{n-1}\left(\left\langle \lambda\bar e_i(0), \bar e_i(0)\right\rangle+
   \int_{0}^{t_{1}}\left(\|\bar e_i'(t)\|^2+\left\langle R\left(\gamma'(t), \bar e_i(t)\right) \gamma'(t), \bar e_i(t)\right\rangle\right)dt\right) \\
&\leq \sum_{i=1}^{n-1}\left(\left\langle \tilde\lambda\hat {\tilde e}_i(0), \hat{\tilde e}_i(0)\right\rangle+
   \int_{0}^{t_{1}}\left(\|\hat {\tilde e}_i'(t)\|^2+\left\langle R\left(\tilde \gamma'(t), \hat {\tilde e}_i(t)\right) \tilde \gamma'(t), \hat {\tilde e}_i(t)\right\rangle\right)dt\right) \\
&=\sum_{i=1}^{n-1}I_{\tilde\lambda\cdot\text{id}}(\hat{\tilde e}_i,\hat{\tilde e}_i)
\end{aligned}
$$
(here $\sum\limits_{i=1}^{n-1}\left\langle R\left(\gamma'(t), \bar e_i(t)\right) \gamma'(t), \bar e_i(t)\right\rangle=-\frac{a^2(t)}{a^2(t_1)}\text{Ric}(\gamma'(t),\gamma'(t))\leq-\frac{a^2(t)}{a^2(t_1)}(n-1)k=\sum\limits_{i=1}^{n-1}\left\langle R\left(\tilde \gamma'(t), \hat {\tilde e}_i(t)\right)\right.$ $\left.\tilde \gamma'(t), \hat {\tilde e}_i(t)\right\rangle$).
I.e., (2.4) holds (by (2.5) and (2.6)); moreover, the equality holds only if $\lambda=\tilde\lambda$ and the Jacobi field $\hat e_i(t)=\bar e_i(t)=\frac{a(t)}{a(t_1)}e_i(t)$ with
$\hat e_i'(0)=\lambda\hat e_i(0)$ on $[0,t_1]$, which implies that $K(\gamma'(t),e_i(t))=k$ and $J_i(t)=a(t_1)\hat e_i(t)=a(t)e_i(t)$ (note that $a(t_1)\hat e_i(0)=e_i(0)=J_i(0)$).

As a result, $\tilde J(t)\geq J(t)$ on $[0,l]$; and if the equality holds for some $t_0\in(0,l]$, then on $[0,t_0]$ each $\|J_i(t)\|=\|\tilde J_{i}(t)\|$, $\frac{J_{i}(t)}{\|J_{i}(t)\|}$ is a parallel vector field, $K(\gamma'(t),J_i(t))=k$,  and $\lambda=\tilde\lambda$.
\hfill$\Box$

\vskip2mm

\noindent {\bf Proof of Corollary E}.

Let $(\rho,\theta_1,\cdots,\theta_{n-1})$ be the polar coordinates of $M_p$. On $M_{\exp_p(\rho,\theta_1,\cdots,\theta_{n-1})}$,
let $J_i(\rho,\theta_1,\cdots,\theta_{n-1})$ denote $\text{d}\exp_p\left(\frac{\partial}{\partial \theta_i}\right)$.
Then we can define a volume form at $(\rho,\theta_1,\cdots,\theta_{n-1})\in M_p$ by
$$\text{d}V=\begin{cases} \|J_1\wedge\cdots\wedge J_{n-1}\|\text{d}\theta_1\cdots \text{d}\theta_{n-1}, & \text{$\exp_p([O(\rho,\theta_1,\cdots,\theta_{n-1})])$ is a minimal geodesic}  \\ 0, & \text{$\exp_p([O(\rho,\theta_1,\cdots,\theta_{n-1})])$ is not a minimal geodesic} \end{cases}$$
(cf. [Pe], [WSY]), where $O$ is the origin of $M_p$; and
$$\text{\rm Vol}(\partial B(p,\rho_0))=\oint \text{d}V|_{\rho=\rho_0},\ \ \text{\rm Vol}(B(p,\rho_0))=\int_{0}^{\rho_0}\oint\text{d}V|_{\rho=r}\text{d}r.$$
Note that each $J_i$ is a Jacobi field along and perpendicular to geodesic $\exp_p([O(\rho,\theta_1,\cdots,\theta_{n-1})])$ (where each $\theta_i$ is fixed). Since $\overline{B(p,r)}$ is isometric to $\overline{B(p',r)}\subset \Bbb S^n_{k'}$, $\|J_1\|_{\rho=r}=\cdots=\|J_{n-1}\|_{\rho=r}$,
and $J_1|_{\rho=r},\cdots,J_{n-1}|_{\rho=r}$ are perpendicular to one another, and $J_i'|_{\rho=r}=\lambda(r) J_i|_{\rho=r}$ (where the derivative is with respect to $\rho$). Similarly,  on $\Bbb S^n_k$, we have Jacobi fields $\tilde J_i(\rho,\theta_1,\cdots,\theta_{n-1})$ along and perpendicular to geodesic $\exp_{\tilde p}([O(\rho,\theta_1,$ $\cdots,\theta_{n-1})])$ so that $\|\tilde J_1\|_{\rho=\tilde r}=\cdots=\|\tilde J_{n-1}\|_{\rho=\tilde r}$, and $\tilde J_1|_{\rho=\tilde r},\cdots,\tilde J_{n-1}|_{\rho=\tilde r}$ are perpendicular to one another, and $\tilde J_i'|_{\rho=\tilde r}=\tilde\lambda(\tilde r) \tilde J_i|_{\rho=\tilde r}$.

Since $\overline{B(p,r)}$ is isometric to $\overline{B(p',r)}\subset \Bbb S^n_{k'}$ with $k'\geq k$ and $\text{\rm Vol}(\partial B(p,r))=\text{\rm Vol}(\partial B(\tilde p,\tilde r))$, we have that $\|J_i\|_{\rho=r}=\|\tilde J_i\|_{\rho=\tilde r}$, and it is not hard to see that $r\geq\tilde r$ with $\tilde r\leq\frac{\pi}{2\sqrt k}$ if $k>0$ and thus $\lambda(r)\leq \tilde \lambda(\tilde r)$ (and equality holds only if $k'=k$, cf. Bishop's volume comparison theorem). Moreover, given any minimal geodesic $[pq]=[p p_0]\cup[p_0q]\subset M$ with $|pp_0|=r$ and $|p_0q|=R$ (note that we can assume $r<\frac{\pi}{\sqrt{k'}}$ if $k'>0$), a point is a focal point of $p_0$ with respect to $\lambda(r)\cdot\text{id}$ on $[p_0q]$ if and only if it is conjugate to $p$, so $q$ is the only possible such a point. Then by taking into account Theorem D we can conclude that, for any $R>0$,
$$\text{\rm Vol}(\partial B(p,r+R))\leq\text{\rm Vol}(\partial B(\tilde p,\tilde r+R)).$$
This implies $\text{\rm Vol}\left(B(p,r+R)\setminus B(p,r)\right)\leq\text{\rm Vol}\left(B(\tilde p,\tilde r+R)\setminus B(\tilde p,\tilde r)\right)$; moreover, if the equality holds,
by Theorem D we have that $\text{Ric}_M=(n-1)k$ and $\lambda(r)=\tilde\lambda(\tilde r)$, and so $k'=k$, $r=\tilde r$,
and $B(p,r+R)$ is isometric to $B(\tilde p,\tilde r+R)$ by the rigidity part of Bishop's volume comparison theorem.
\hfill$\Box$

\vskip2mm

\noindent {\bf Remark 2.1.} From the proof of Corollary E, it is not hard to see that the condition `$\text{\rm Vol}(\partial B(p,r))=\text{\rm Vol}(\partial B(\tilde p,\tilde r))$' in the corollary can be replaced with `$\text{\rm Vol}(\partial B(p,r))\leq\text{\rm Vol}(\partial B(\tilde p,\tilde r))$ with $r\geq\tilde r$'.

\vskip8mm

\noindent School of Mathematical Sciences (and Lab. Math. Com.
Sys.), Beijing Normal University, Beijing, 100875,
People's Republic of China

\noindent E-mail: 13716456647@163.com; suxiaole$@$bnu.edu.cn; wyusheng$@$bnu.edu.cn

\end{document}